\documentclass[psamsfonts,reqno]{amsart}

\usepackage{amssymb}
\usepackage{amscd}
\usepackage{eucal}
\usepackage{psfig}
\usepackage[dvips]{graphicx}
\usepackage[small]{caption}

\numberwithin{equation}{section}

\DeclareMathOperator{\lm}{\text{l.i.m.}}

\newcommand{\beq}{\begin{equation}}
\newcommand{\eeq}{\end{equation}}

\newcommand{\Real}{\mbox{Re\,}}
\newcommand{\Imag}{\mbox{Im\,}}

\newcommand{\R}{\mathbb{R}}
\newcommand{\Z}{\mathbb{Z}}
\newcommand{\C}{\mathbb{C}}

\newcommand{\HH}{\mathbb{H}}

\newcommand{\cA}{\mathcal{A}}

\newcommand{\cI}{\mathcal{I}}

\newcommand{\cP}{\mathcal{P}}

\newcommand{\uf}{\underline{f}}
\newcommand{\Cons}{\text{\rm const.}}

\newcommand{\sfrac}[2]{{\vphantom1\smash{\lower.5ex\hbox{\small$#1$}}\over
        \vphantom1\smash{\raise.4ex\hbox{\small$#2$}}}} 

\newtheorem*{theorem-nonumber}{Theorem}
\newtheorem{proposition}{Proposition}

\theoremstyle{remark}
\newtheorem*{remark}{Remark}
\newtheorem*{remarks}{Remarks}

\theoremstyle{definition}

\begin{document}

\title[Holomorphic Extension]{Holomorphic Extension associated with Fourier--Legendre Expansions}

\author[E. De Micheli]{E. ~De Micheli}
\address[E. De Micheli]{IBF -- Consiglio Nazionale delle Ricerche \\ Via De Marini, 6 - 16149 Genova, Italy}
\email[E.~De Micheli]{demicheli@ge.cnr.it}

\author[G. A. Viano]{G. A. ~Viano}
\address[G. A. ~Viano]{Dipartimento di Fisica - Universit\`a di Genova,
Istituto Nazionale di Fisica Nucleare - sez. di Genova\\
Via Dodecaneso, 33 - 16146 Genova, Italy}
\email[G.A.~Viano]{viano@ge.infn.it}

\begin{abstract}
In this article we prove that if the coefficients of a
Fourier--Legendre expansion satisfy a suitable Hausdorff--type
condition, then the series converges to a function which admits a
holomorphic extension to a cut--plane. Furthermore, we prove that
a Laplace--type (Laplace composed with Radon) transform of the
function describing the jump across the cut is the unique
Carlsonian interpolation of the Fourier coefficients of the
expansion. We can thus reconstruct the discontinuity function from
the coefficients of the Fourier--Legendre series by the use of the
Pollaczek polynomials.
\end{abstract}

\maketitle

\section{Introduction}
\label{se:introduction}
Let us consider the following Fourier--Legendre series
\beq
\label{I1}
\frac{1}{4\pi}\sum_{n=0}^\infty (2n+1) a_n P_n(\cos\theta).
\eeq
The classical theory of polynomial expansions, as described by Walsh \cite{Walsh},
establishes for these expansions convergence properties which are closely analogous
to the well--known convergence properties of the Taylor series expansions. In this
case the region of convergence, instead of being circles of radius $\rho$, are
ellipses $E_\rho$ with foci $\pm 1$ and \emph{radius} $\rho=$(semiminor+semimajor)--axis.
The expansion converges inside the largest ellipse within which the function being
expanded in terms of the series (\ref{I1}) is holomorphic.

The following question quite naturally arises: Is it possible to find suitable
conditions on the coefficients $a_n$ which allow a holomorphic extension of the
function (to which the series (\ref{I1}) converges) to the whole complex
$\cos\theta$--plane except for a cut along the positive axis? The answer to this
question is positive and quite analogous to that derived in \cite{DeMicheli1}
in connection with the Taylor and Laurent series: essentially, the coefficients
$\{a_n\}$ are required to satisfy suitable Hausdorff--type conditions.

\begin{figure}[tb]
\begin{center}
\includegraphics[width=9cm]{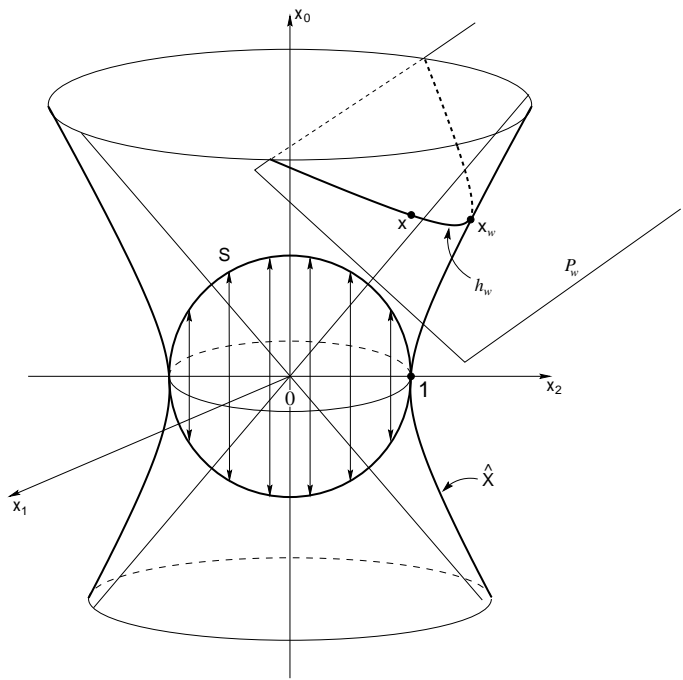}
\caption{\label{figura_1}Horocyclic fibration of the one--sheeted hyperboloid.}
\end{center}
\end{figure}

To prove these results it is convenient to proceed through two steps. First,
replacing the complex $\cos\theta$--plane ($\theta\in\C;\,\theta=u+iv;\,u,v\in\R$)
by a complex hyperboloid $X^{(c)}$, which contains as submanifolds the
\emph{Euclidean sphere} $S=(i\R\times\R^2)\cap X^{(c)}$ which gives the support of
the SO(3,$\R$) harmonic analysis, and the real one--sheeted hyperboloid
$X=\R^3\cap X^{(c)}$ that contains the support of the cut (see Fig. \ref{figura_1}).
In the second step we consider a fibration on a meridian hyperbola of $X^{(c)}$, which
is obtained through a Radon--type transformation. This fibration allows us to reduce
the harmonic analysis to that associated with a complex one--dimensional hyperbola,
which \emph{contains} the Euclidean--circle and the real hyperbola. We are thus led
to regard the series (\ref{I1}) as a trigonometrical series on the Euclidean--circle,
making it possible to apply several results of the same type as those obtained in
\cite{DeMicheli1}. In particular, we can prove that this trigonometrical series
converges to a function which admits a holomorphic extension to a complex cut--plane
if the coefficients $a_n$ satisfy a suitable Hausdorff--type condition. Then, by
inverting the Radon transform, we return to the complex $\theta$--plane and,
finally, study the holomorphic extension associated with the Legendre series (\ref{I1}).
In addition, we obtain a unique Carlsonian interpolation of the $a_n$'s, denoted
by $\tilde{a}(\lambda)$, which turns out to be the composition of the ordinary Laplace
transform with the Radon transform. By inverting these Laplace and Radon transforms,
we are then able to reconstruct the jump function across the cut using the Pollaczek
polynomials.

In a classical article \cite{Stein}, Stein and Wainger have already proved the
holomorphic extension associated with the series (\ref{I1}) (see, in particular,
Theorem 4 in their article). However, our work differs from Stein and Wainger's for
several reasons:
\begin{itemize}
\item[i)] We associate the holomorphic extension of the Legendre series to the
Hausdorff condition on the coefficients $a_n$, and, accordingly, we introduce the
Carlsonian interpolation $\tilde{a}(\lambda)$ of these coefficients which is the
composition of an ordinary Laplace transform with a Radon transform (called
\emph{spherical--Laplace transform} in \cite{Faraut}).
\item[ii)] The inversion of this \emph{Laplace composed with Radon} transform allows
us to reconstruct the jump function across the cut starting from the coefficients
$\{a_n\}$ of the series (\ref{I1}). This procedure is extremely relevant in the
inverse problem in quantum scattering theory for the class of Yukawian potentials
(see \cite{Fioravanti}). In fact, from the discontinuity function across the cut
the spectral density associated with the Yukawian potentials can be determined.
\item[iii)] For particular values of $\lambda$, i.e., $\lambda=-1/2+i\nu$ $(\nu\in\R)$,
we obtain from $\tilde{a}(\lambda)$ the classical Mehler transform, which is
precisely the mathematical tool used by Stein and Wainger for obtaining their results,
through the Plancherel theorem.
\item[iv)] Finally, we can give a geometrical interpretation of our results and methods
by introducing the complex hyperboloid $X^{(c)}$ and associating to it, through a Radon
transform, a fibration on a meridian hyperbola.
\end{itemize}

This paper is organized as follows: in Section \ref{se:interpolation} we study
the Carlsonian interpolation of the Hausdorff moments and, correspondingly,
the Hardy spaces to which this interpolation belongs. In Section \ref{se:legendre}
we prove that the Legendre expansion can be regarded as a trigonometrical series.
In Section \ref{se:holomorphic} we prove a holomorphic extension associated with
these trigonometrical series in the complex $\tau$--plane using the same procedure
adopted in \cite{DeMicheli1}. The variable $\tau$ can then be interpreted as one
of the horocyclic coordinates related to the fibration on the meridian hyperbola
$\hat{X}^{(c)}$ (see the Appendix). In Section \ref{se:inversion} we study the
inversion of the Radon--Abel transformation and we prove the holomorphic extension
in the $\theta$--plane associated with the Legendre series. In Section
\ref{se:laplace} we show that the Carlsonian interpolation $\tilde{a}(\lambda)$
of the Hausdorff moments can be represented as the composition of the ordinary
Laplace transform with the Radon transform, and we find an integral representation
of the jump function across the cut; moreover, we show how to reconstruct the
discontinuity function (across the cut) starting from the Fourier--Legendre
coefficients $\{a_n\}$, using the Pollaczek polynomials. Finally, in the Appendix
we illustrate all the geometrical aspects of the method we used.

\section{Interpolation of Hausdorff moments and Hardy spaces}
\label{se:interpolation} Let us consider a sequence
$\{f_n\}_0^\infty$ of (real) numbers $f_n$, and denote by $\Delta$
the difference operator
\beq
\label{uno}
\Delta f_n = f_{n+1} - f_n.
\eeq
Then we have:
\beq
\label{due}
\Delta^k f_n=
\underbrace{\Delta\times\Delta\times\cdots\times\Delta}_{k} f_n =
\sum_{m=0}^k (-1)^m {k\choose m} f_{n+k-m},
\eeq
(for every $k \geq 0$); $\Delta^0$ is the identity operator by definition. Now,
suppose that there exists a positive constant $M$ such that:
\beq
\label{tre}
(n+1)^{(1+\epsilon)}\sum_{i=0}^n {n\choose i}^{(2+\epsilon)} \left | \Delta^i f_{(n-i)} \right|^{(2+\epsilon)} < M
\qquad (n=0,1,2,\ldots; \epsilon>0).
\eeq
It can be proved \cite{Widder} that condition (\ref{tre}) is
necessary and sufficient to represent the sequence
$\left\{f_n\right\}_0^\infty$ as follows:
\beq
\label{quattro}
f_n = \int_0^1 x^n u(x) \, dx \qquad (n=0,1,2,\ldots),
\eeq
where $u(x)$ belongs to $L^{2+\epsilon}[0,1]$.

We can prove the following proposition.

\begin{proposition}
\label{pro1}
Let the sequence $\left\{f_n\right\}_0^\infty$, $f_n=n^p a_n$, $(p\geq 1)$,
satisfy condition $(\ref{tre})$. Then there exists a unique Carlsonian
interpolation of the numbers $a_n$, denoted by $\tilde{a}(\lambda)$
$(\lambda \in \C, [\tilde{a}(\lambda)]_{(\lambda=n)} = a_n, n=0,1,2,\ldots)$,
that satisfies the following conditions:
\begin{itemize}
\item[i)] $\tilde{a}(\lambda)$ is holomorphic in the half--plane $\Real\lambda > -1/2$,
continuous at $\Real\lambda = -1/2$;
\item[ii)] $\lambda^p \tilde{a}(\lambda)$ belongs to $L^2(-\infty,+\infty)$ for
any fixed value of $\Real\lambda \geq -1/2$: i.e., putting $\lambda=\sigma+i\nu$,
\beq
\label{new1}
\int_{-\infty}^{+\infty}\left|(\sigma+i\nu)^p\,\tilde{a}(\sigma+i\nu)\right|^2\,d\nu < \infty;
\eeq
\item[iii)] $\lambda^p \tilde{a}(\lambda)$ tends uniformly to zero as $\lambda$
tends to infinity inside any fixed half--plane $\Real\lambda \geq \delta > -1/2$;
\item[iv)] $\lambda^{(p-1)} \tilde{a}(\lambda)$ belongs to $L^1(-\infty,+\infty)$
for any fixed value of $\Real\lambda \geq -1/2$.
\end{itemize}
\end{proposition}

\begin{proof}
If the sequence $\{f_n\}_0^\infty$ satisfies condition (\ref{tre}),
then representation (\ref{quattro}) holds true. If we put $x=e^{-t}$ in the integral
of (\ref{quattro}) we obtain:
\beq
\label{cinque}
f_n = \int_0^{+\infty} e^{-nt} e^{-t} u(e^{-t}) \, dt \qquad (n=0,1,2,\ldots).
\eeq
Therefore the numbers $f_n$ can be regarded as the restriction to the integers of the
following Laplace transform:
\beq
\label{sei}
\tilde{F}(\lambda) = \int_0^{+\infty} e^{-(\lambda+1/2)t} e^{-t/2} u(e^{-t}) \, dt.
\eeq
It can easily be verified that $[\tilde{F}(\lambda)]_{(\lambda=n)}=f_n$.
By applying the Paley--Wiener theorem to equality (\ref{sei}), and
recalling that the function $\exp(-t/2)u(\exp(-t))$ belongs to $L^2[0,+\infty)$, we can
conclude that $\tilde{F}(\lambda)$ belongs to the Hardy space $\HH^2(\C_{-1/2})$
$(\C_{-1/2}=\{\lambda\in\C,\Real\lambda>-1/2\})$. We can thus apply the Carlson
theorem \cite{Boas}, and state that $\tilde{F}(\lambda)$ is the unique Carlsonian
interpolation of the numbers $f_n$.
Furthermore, by noting that $\tilde{F}(\lambda)=\lambda^p\tilde{a}(\lambda)$ $(p\geq 1)$,
properties (ii), (iii) and the analyticity of $\tilde{a}(\lambda)$ in the half--plane
$\Real\lambda > -1/2$ follow.
Next let us note that the function $\exp(-t/2)u(\exp(-t))$ belongs to $L^1[0,+\infty)$;
in fact, we can state that
$\int_0^{+\infty}|\exp(-t/2)u(\exp(-t))|\,dt=\int_0^1|u(x)/\sqrt{x}|\,dx<\infty$,
in view of the fact that $u \in L^{2+\epsilon}[0,1]$.
Therefore, from the Riemann--Lebesgue theorem applied to representation (\ref{sei})
it follows that the function $\tilde{F}(-1/2+i\nu)$ $(\nu\in\R)$ is continuous,
and thus property (i) is proved.\\
Concerning property (iv) we may use the Schwarz inequality and write
\beq
\label{seiprimo}
\begin{split}
&\int_{-\infty}^{+\infty}\left |(\sigma+i\nu)^{(p-1)} \tilde{a}(\sigma+i\nu) \right |\,d\nu =
\int_{-\infty}^{+\infty} \left |\frac{\tilde{F}(\sigma+i\nu)}{(\sigma+i\nu)}\right |\,d\nu \\
& \qquad\leq\left ( \int_{-\infty}^{+\infty} \frac{1}{|(\sigma+i\nu)|^2}\,d\nu\right )^{1/2}
\left ( \int_{-\infty}^{+\infty} |\tilde{F}(\sigma+i\nu)|^2\, d\nu \right )^{1/2} < \infty,
\end{split}
\eeq
if $\sigma\geq -1/2$, $\sigma\neq 0$, $p\geq 1$. In fact, let us note that
$\tilde{F}(\sigma+i\nu)\in L^2(-\infty,+\infty)$ for any fixed value of
$\Real\lambda=\sigma\geq -1/2$.
Finally, in view of the regularity and integrability of the function
$\lambda^{(p-1)}\tilde{a}(\lambda)$ in the neighborhood of $\Real\lambda=0$,
we can conclude that $\lambda^{(p-1)}\tilde{a}(\lambda)$ belongs to $L^1(-\infty,+\infty)$
for any fixed value of $\Real\lambda=\sigma\geq-1/2$, $(p\geq 1)$.
\end{proof}

\begin{remark}
In order to prove the continuity of $\tilde{a}(\lambda)$ at
$\Real\lambda=-1/2$ ($\lambda=-1/2+i\nu,\,\nu\in\R$), it is necessary to use
condition (\ref{tre}) which is slightly more restrictive than condition (8) of
\cite{DeMicheli1} where the term $\epsilon$ was missing.
\end{remark}

\section{Legendre expansions as trigonometrical series}
\label{se:legendre}
Let us consider the following Legendre series:
\beq
\label{sette}
\frac{1}{4\pi}\sum_{n=0}^\infty (2n+1)\,a_n\,P_n(\cos u),
\eeq
where $P_n$ denotes the Legendre polynomials. The polynomials $P_n$ satisfy the
following integral representation:
\beq
\label{otto}
P_n(\cos u) = \frac{1}{\pi}\int_0^\pi (\cos u+i\sin u\cos\eta)^n\,d\eta.
\eeq
Let us now suppose that expansion (\ref{sette}) converges to a function
$\uf(\cos u)$ but, for the moment, we shall leave the type of convergence unspecified.
We only assume that $\uf(\cos u)$ is a measurable and integrable function in the
interval $u\in[0,\pi]$. Thus, the Legendre coefficients $a_n$ can be written as
\beq
\label{nove}
a_n = 2\pi\int_0^\pi \uf (\cos u) P_n(\cos u)\sin u\,du.
\eeq
Our goal now is to rewrite expansion (\ref{sette}) as a trigonometrical series.
For this purpose, we prove the following proposition.

\begin{proposition}
\label{pro2}
The Legendre coefficients $\{a_n\}_0^\infty$ coincide with the Fourier coefficients
of the form:
\beq
\label{dieci}
a_n = \int_{-\pi}^\pi \hat{f}(t) e^{int}\,dt \qquad (n=0,1,2,\ldots),
\eeq
where
\beq
\label{undici}
\hat{f}(t) = -2i\epsilon(t)e^{it/2}\int_0^t f(u)\left[2(\cos u-\cos t)\right]^{-1/2}\sin u\,du,
\eeq
with $f(u)\equiv\uf(\cos u)$, and $\epsilon(t)$ being the sign function.
\end{proposition}

\begin{proof}
From the Dirichlet--Murphy integral representation of the Legendre
polynomials (see Ref. \cite{Vilenkin}, Ch. III, Section 5.4):
\beq
\label{diciotto}
P_n(\cos u)=-\frac{i}{\pi}\int_u^{(2\pi-u)}e^{i(n+1/2)t}\left[2(\cos u - \cos t)\right]^{-1/2}\,dt,
\eeq
and from equality (\ref{nove}), we have
\beq
\label{diciannove}
\frac{ia_n}{2}=
\int_0^\pi du\,\uf(\cos u)\sin u\int_u^{(2\pi-u)}e^{i(n+1/2)t}
\left[2(\cos u - \cos t)\right]^{-1/2}\,dt.
\eeq
Inverting the order of integration in formula (\ref{diciannove}), we get
\beq
\label{venti}
\begin{split}
\frac{ia_n}{2}&=\int_0^\pi dt \, e^{i(n+1/2)t}
\int_0^t du\,\sin u\,\uf(\cos u)\left[2(\cos u - \cos t)\right]^{-1/2} \\
&\qquad +\int_\pi^{2\pi} dt\, e^{i(n+1/2)t}
\int_0^{(2\pi-t)} du\,\sin u\,\uf(\cos u) \left[2(\cos u - \cos t)\right]^{-1/2}.
\end{split}
\eeq
Next, from the second integral in the r.h.s. of formula (\ref{venti}), if we
perform the following change of variable, $t\rightarrow t-2\pi$, and change
$u \rightarrow -u$, we get
\begin{equation*}
e^{i\pi}\int_{-\pi}^0 dt\, e^{i(n+1/2)t}\int_0^t du\,\sin u\,\uf(\cos u)
\left[2(\cos u - \cos t)\right]^{-1/2}.
\end{equation*}
Finally, we obtain
\beq
\label{ventuno}
\begin{split}
\frac{ia_n}{2}&=\int_0^\pi dt \, e^{i(n+1/2)t}
\int_0^t du\,\sin u\,\uf(\cos u)\left[2(\cos u - \cos t)\right]^{-1/2} \\
&\qquad +e^{i\pi}\int_{-\pi}^0 dt\, e^{i(n+1/2)t}\int_0^t du\,\sin u\,\uf(\cos u)
\left[2(\cos u - \cos t)\right]^{-1/2},
\end{split}
\eeq
which gives
\beq
\label{ventidue}
a_n=\int_{-\pi}^\pi \hat{f}(t) \,e^{int}\,dt,
\eeq
with $\hat{f}(t)$ given by (\ref{undici}).\\
(For a proof of this result in a more general setting see also \cite{Bros}III).
\end{proof}

It can easily be verified that (see formula (\ref{undici}))
\beq
\label{ventiquattro}
\hat{f}(t)= -e^{it}\hat{f}(-t),
\eeq
and, accordingly, from (\ref{ventidue}) and (\ref{ventiquattro}) we have
\beq
\label{venticinque}
a_n = -a_{-n-1} \qquad (n\in\Z).
\eeq
We are thus prompted to consider the following trigonometrical series,
\beq
\label{ventisei}
\begin{split}
\frac{1}{2\pi}\sum_{n=-\infty}^{+\infty} a_n\,e^{-int}&=
\frac{1}{2\pi}\left\{\sum_{n=0}^{+\infty} a_n\,e^{-int}-e^{it}
\sum_{n=0}^{+\infty} a_n\,e^{int}\right\}  \\
&=\frac{1}{2\pi}\,e^{i\frac{(t-\pi)}{2}}
\sum_{n=-\infty}^{+\infty}(-1)^n a_n\,\cos\left[\left(n+\frac{1}{2}\right)(t-\pi)\right] \\
&= \frac{1}{2\pi}\,e^{i\frac{(t-\pi)}{2}}
\sum_{n=-\infty}^{+\infty} a_n\,\sin\left[\left(n+\frac{1}{2}\right)t\right],
\end{split}
\eeq
and study the holomorphic extension associated with it.

\section{Holomorphic extension associated with the trigonometrical series}
\label{se:holomorphic}
In the complex plane $\C$ of the variable $\tau=t+iw$ $(t,w\in\R)$ we introduce
the following domains:
${\cI}_+^{(\pm\xi_0)} = \{\tau \in \C \mid \Imag\tau > \pm \xi_0,\,\xi_0 \geq 0\}$,
and
${\cI}_-^{(\pm\xi_0)} = \{\tau \in \C \mid \Imag\tau < \pm \xi_0,\,\xi_0 \geq 0\}$.
Correspondingly, we introduce the following cut--domains:
${\cI}_+^{(\xi_0)}\setminus \Xi_+^{(\xi_0)}$, where
$\Xi_+^{(\xi_0)}=\{\tau\in\C \mid\tau=2k\pi + iw,\, w>\xi_0,\, \xi_0\geq 0,\, k\in \Z\}$,
and
${\cI}_-^{(\xi_0)}\setminus \Xi_-^{(-\xi_0)}$, where
$\Xi_-^{(-\xi_0)}=\{\tau\in\C \mid\tau=2k\pi + iw,\, w<-\xi_0,\, \xi_0\geq 0,\, k\in \Z\}$.
We shall use the notation $\dot{A}=A\setminus 2\pi\Z$ for every subset $A$ of $\C$
which is invariant under the translation group $2\pi\Z$.
We can then prove the following proposition.

\begin{proposition}
\label{pro3}
Let us consider the following trigonometrical series,
\beq
\label{ventisette}
\frac{1}{2\pi}\sum_{n=0}^\infty a_n \, e^{-in\tau} \qquad (\tau=t+iw;\,t,w\in\R),
\eeq
and suppose that the set of numbers $\{f_n\}_0^\infty$, $f_n=n^p a_n$,
$(n=0,1,2,\ldots, p\geq 1)$ satisfies condition $(\ref{tre})$. Then:
\begin{itemize}
\item[i)] The series $(\ref{ventisette})$ converges to a function
$\hat{f}(\tau)$ holomorphic in $\cI_-^{(0)}$, the convergence being uniform in
any compact subdomain of $\cI_-^{(0)}$.
\item[ii)] The function $\hat{f}(\tau)$ admits a holomorphic extension to the
cut--domain ${\cI}_+^{(0)}\setminus \dot{\Xi}_+^{(0)}$, i.e., it is analytic in
$\C \setminus \{\tau=2k\pi+iw\,|\,k\in\Z,w >0\}$.
\item[iii)] The jump function $\hat{F}(w)$ (which equals the discontinuity of
$\,i\hat{f}(\tau)$ across the cuts $\dot{\Xi}_+^{(0)}$) is a function of class
$C^{p-1}$, $(p\geq 1)$, and satisfies the following bound
\beq
\label{ventotto}
\left|\hat{F}(w)\right|\leq\left\|\tilde{a}_\sigma\right\|_1\,e^{\sigma w},
\qquad\left(\sigma\geq-\frac{1}{2}, w\in\R^+\right),
\eeq
where $\tilde{a}(\sigma+i\nu)$ $(\nu\in\R)$ is the Carlsonian interpolation of the
coefficients $a_n$, and
\beq
\label{ventinove}
\left\|\tilde{a}_\sigma\right\|_1=
\frac{1}{2\pi}\int_{-\infty}^{+\infty}|\tilde{a}(\sigma+i\nu)|\,d\nu,
\qquad\left(\sigma\geq-\frac{1}{2}\right).
\eeq
\item[iv)] $\tilde{a}(\sigma+i\nu)$ is the Laplace transform of the jump function
$\hat{F}(w)$: i.e.,
\beq
\label{trenta}
\tilde{a}(\sigma+i\nu)=
\int_0^{+\infty}\hat{F}(w)\,e^{-(\sigma+i\nu)w}\,dw, \qquad\left(\sigma>-\frac{1}{2}\right).
\eeq
\item[v)] The following Plancherel equality holds true:
\beq
\label{trentuno}
\int_{-\infty}^{+\infty}\left |\tilde{a}(\sigma+i\nu)\right |^2\,d\nu=
2\pi\int_{-\infty}^{+\infty}\left |\hat{F}(w)\,e^{-\sigma w}\right |^2\,dw,
\qquad\left(\sigma\geq-\frac{1}{2}\right).
\eeq
\end{itemize}
\end{proposition}

\begin{proof}
See the proof of Theorem 1 in \cite{DeMicheli1}
(see also the remark at the end of Proposition \ref{pro1}).
\end{proof}

\begin{remark}
Let us note that the Plancherel equality (\ref{trentuno}) holds true
under the milder condition that the coefficients $a_n$ ($n=0,1,2,\ldots$) satisfy
condition (\ref{tre}).
\end{remark}

Next, we can state the following proposition.

\begin{proposition}
\label{pro4} If in the trigonometrical series
\beq
\label{trentadue}
\frac{1}{2\pi}\left\{\sum_{n=0}^\infty
a_n\,e^{-int}-e^{it} \sum_{n=0}^\infty a_n\,e^{int}\right\} \qquad (t\in\R),
\eeq
the coefficients $a_n$ satisfy the assumptions
required by Proposition $\ref{pro3}$, then:
\begin{itemize}
\item[i)] the series converges to a continuous function $\hat{f}(t)$, $(t\in\R)$,
the convergence being uniform on any compact subdomain of the real line.
\item[ii)] The function $\hat{f}(t)$ admits a holomorphic extension to the cut--domain
${\cI}_+^{(0)}\setminus\dot{\Xi}_+^{(0)} \cup {\cI}_-^{(0)}\setminus\dot{\Xi}_-^{(0)}$:
i.e., it is analytic in $\C \setminus \{\tau = 2k\pi+iw \,|\, k \in\Z, |w|>0\}$.
\item[iii)] The jump function across the cuts $\dot{\Xi}_\pm^{(0)}$ satisfies
conditions analogous to (iii)--(v) of Proposition $\ref{pro3}$.
\end{itemize}
\end{proposition}

\begin{proof}
Statement (i) follows from observing that, in view of the assumptions
on the coefficients $a_n$ we have:
\beq
\label{trentatre}
\left |\frac{1}{2\pi}\left\{\sum_{n=0}^\infty a_n\,e^{-int}-e^{it}
\sum_{n=0}^\infty a_n\,e^{int}\right\}\right | \leq \frac{1}{\pi}\sum_{n=0}^\infty |a_n|<\infty.
\eeq
Next, by the Weierstrass theorem on the uniformly convergent series of continuous
functions, we obtain the result. Statements (ii) and (iii) can be proved analogously
to the proof of the corresponding statements in Proposition \ref{pro3}.
\end{proof}

\begin{remarks}
(i) If the coefficients $a_n$ in (\ref{trentadue}) are exponentially
bounded, i.e., $|a_n|\leq K\exp(-(n-m)\xi_0)$, ($n>m,\,m\in\R^+,\,\xi_0>0,\,K=\,{\rm constant}$),
then $\hat{f}(t)$ admits a holomorphic extension to the cut--domain
${\cI}_+^{(0)}\setminus \dot{\Xi}_+^{(\xi_0)} \cup {\cI}_-^{(0)}\setminus \dot{\Xi}_-^{(-\xi_0)}$
(see Proposition 5 of \cite{DeMicheli1}).
For the sake of simplicity, in the following we shall only consider the case $\xi_0=0$. \\
(ii) Similarly, we assume hereafter that condition (\ref{tre}) is satisfied by the
whole sequence $\{f_n\}_0^\infty$; we could also assume that this condition is
satisfied only by the subset $\{f_n\}_{n_0}^\infty$, $(n_0>0)$, $f_n=n^p a_n$, $p\geq 1$.
In this case the result proved above still holds but for minor modifications,
e.g., that in (\ref{ventotto}) now $\sigma\geq(n_0-\frac{1}{2})$, $(n_0>0)$,
and likewise in formulae (\ref{ventinove}), (\ref{trenta}), and
(\ref{trentuno}). See also the remark after Proposition 5 in \cite{DeMicheli1}.
\end{remarks}

\section{Inversion of the Radon--Abel transformation
and holomorphic extension associated with the Legendre Series}
\label{se:inversion}
Proposition \ref{pro2} allows us to regard the Legendre expansions as trigonometrical
series. The function $\hat{f}(t)$ is the Radon--Abel transformation of the function
$f(u)$ and, moreover, it can be regarded as the restriction of a function $\hat{f}(\tau)$,
$(\tau\in\C)$, which is the Radon--Abel transformation of a function $f(\theta)$,
$(\theta\in\C)$, when $\tau=t$ and $\theta=u$ ($t,\,u\in\R$) (see the Appendix).
It is therefore of primary interest to derive the inversion of the Radon--Abel transformation.
We can prove the following proposition.

\begin{proposition}
\label{pro5}
Let us suppose that the sequence $f_n=n^p a_n$, $(p\geq 2)$, (the numbers $a_n$
being the coefficients of the Legendre expansion $(\ref{sette}))$ satisfies
the Hausdorff condition $(\ref{tre})$. We can then write the following Radon--Abel
transformation (see formula $(\ref{a.20}$), and $(\ref{a.21})$ of the Appendix),
\beq
\label{trentaquattro}
\hat{f}(t)=-2e^{it/2}\int_0^t f(u)\left[2(\cos t-\cos u)\right]^{-1/2}\sin u\,du,
\eeq
which admits the following inversion:
\beq
\label{trentacinque}
f(u)=\frac{1}{\pi\sin u}\frac{d}{du}
\int_0^u e^{-it/2}\hat{f}(t)\left[2(\cos u-\cos t)\right]^{-1/2}\sin t\,dt.
\eeq
\end{proposition}

\begin{proof}
In view of the assumptions on the Legendre coefficients $a_n$ and
of Propositions \ref{pro2}, \ref{pro3} and \ref{pro4} we can guarantee that
representation (\ref{trentaquattro}) holds true, and, moreover, that the function
$\hat{f}(t)$, $(t\in\R)$, is continuous. Furthermore, since the sequence
$n^p a_n$ satisfies the Hausdorff condition (\ref{tre}) with $p\geq 2$, $\hat{f}(t)$
is also differentiable. This fact can easily be proved by differentiating the series
at the l.h.s. of formula (\ref{ventisei}) term by term, and observing that it is
majorized by the convergent series: $\Cons\sum_{n=0}^\infty |n a_n|<\infty$.\\
Let us note that, for the sake of simplicity, we work with representation (\ref{trentaquattro})
instead of (\ref{undici}) (see also the Appendix, formulae (\ref{a.20}) and (\ref{a.21})).

Next, we set in formula (\ref{trentaquattro}) $\cos t=(1-\rho)$, ($\rho>0$),
$\cos u=(1-\rho')$, ($0\leq\rho'\leq\rho$). Thus formula (\ref{trentaquattro}) can be
rewritten as
\beq
\label{trentasei}
\hat{f}(t)=2ie^{it/2}\int_0^\rho \uf(1-\rho')\left[2(\rho-\rho')\right]^{-1/2}\,d\rho'.
\eeq
Then we introduce the Riemann--Liouville integral $[I_\alpha\phi]$, which can be
written as
\beq
\label{trentasette}
[I_\alpha\phi](\rho)=
\frac{1}{\Gamma(\alpha)}\int_0^\rho\phi(\rho')(\rho-\rho')^{(\alpha-1)}\,d\rho' \qquad (\alpha>0),
\eeq
and we have
\beq
\label{trentotto}
\hat{f}(t)=i\sqrt{2\pi}e^{it/2} [I_{1/2}\phi](t),
\eeq
where
\beq
\label{trentanove}
\phi(\rho')=\uf(1-\rho').
\eeq
If $[I_\alpha\phi]$ is $\alpha$--times differentiable, then the following properties
of the Riemann--Liouville integral can be applied:
\begin{eqnarray}
\label{quaranta}
\mbox{i)} &~& \quad I_\alpha \circ I_\beta = I_{\alpha+\beta} \quad (\alpha,\beta > 0) \quad
(\text{\rm in particular,}~ I_{1/2} \circ I_{1/2} = I_1).\label{quaranta} \\
\label{quarantuno}
\mbox{ii)} &~& \quad \left(\frac{d}{d\rho}\right)^\alpha [I_\alpha\phi](\rho) = \phi(\rho)
\qquad (\alpha=1,2,\ldots).
\end{eqnarray}
Since $\hat{f}(t)$ is differentiable we can write, by using properties (i) and (ii):
\beq
\label{quarantadue}
\frac{d}{d\rho'}\left[I_{1/2}[I_{1/2}\phi]\right] = \frac{1}{\sqrt\pi}\frac{d}{d\rho'}
\int_0^{\rho'} [I_{1/2}\phi](\rho)(\rho'-\rho)^{-1/2}d\rho = \phi(\rho').
\eeq
By applying the last equalities to our case, and in view of (\ref{trentotto}), we
obtain formula (\ref{trentacinque}).
\end{proof}

For what concerns the inversion of the Radon--Abel transformation in a more general
setting see Ref. \cite{Abouelaz}.

Proposition \ref{pro4} proves that, if the numbers $n^p a_n$ ($p\geq 1$) satisfy
the Hausdorff condition (\ref{tre}), then $\hat{f}(t)$ is the restriction to the
real axis of a function $\hat{f}(\tau)$ holomorphic in
$\dot{{\cI}}= \left({\cI}_+^{(0)}\setminus \dot{\Xi}_+^{(0)}\right) \cup
\left({\cI}_-^{(0)}\setminus \dot{\Xi}_-^{(0)}\right)$.
Assuming hereafter that the Hausdorff condition (\ref{tre}) is satisfied by the sequence
$f_n=n^p a_n$ with $p\geq 2$, we can extend representation (\ref{trentacinque})
uniquely in the following way:
\beq
\label{quarantaquattro}
f(\theta)=\frac{1}{\pi\sin\theta}\frac{d}{d\theta}\int_{\gamma_\theta}
e^{-i\tau/2}\hat{f}(\tau)\frac{\sin\tau}{[2(\cos\theta-\cos\tau)]^{1/2}}\,d\tau,
\eeq
where $\gamma_\theta$ denotes the ray $\underline{\gamma}_\theta$ oriented
from 0 to $\theta$. We can now prove the following proposition.

\begin{proposition}
\label{pro6}
Let us suppose that the sequence $f_n=n^p a_n$ $(n=0,1,\ldots)$ satisfies the Hausdorff
condition $(\ref{tre})$ with $p\geq 2$, then the function $f(\theta)$ represented
by formula $(\ref{quarantaquattro})$ is even, $2\pi$--periodic, and holomorphic
in $\dot{{\cI}}$ (here referred to the complex plane of the variable $\theta=u+iv$).
\end{proposition}

\begin{proof}
The assumptions on the Legendre coefficients $a_n$ allow us to state
that $\hat{f}(\tau)$ is a $2\pi$--periodic function holomorphic in the domain
$\dot{{\cI}}$ of the complex $\tau$--plane (see Proposition \ref{pro4}) that also
satisfies the following symmetry property:
\beq
\label{quarantacinque}
\hat{f}(\tau)=-e^{i\tau}\,\hat{f}(-\tau),
\eeq
which derives from equality (\ref{ventiquattro}) in view of the uniqueness of the
analytic continuation (see also the Appendix, formula (\ref{a.17})).
The above--mentioned properties imply that $\hat{f}(\tau)$ is of the following form:
$\hat{f}(\tau)=\exp(i\tau/2)(1-\cos\tau)^{1/2}\,b(\cos\tau)$, with $b(\cos\tau)$ analytic in
$\underline{D}=\{\cos\tau\in\C,\tau\in\dot{{\cI}}\}$.
Through the following parametrization of $\gamma_\theta$: $\cos\tau=1+\lambda(\cos\theta-1)$,
$(0\leq\lambda\leq 1)$ (see also the Appendix), the r.h.s. of formula (\ref{quarantaquattro})
can be rewritten as
\beq
\label{quarantasei}
\frac{i}{\sqrt{2}\pi}\frac{d}{d(\cos\theta)}\left\{(\cos\theta -1)
\int_0^1 b[1+\lambda(\cos\theta-1)]\lambda^{1/2}(1-\lambda)^{-1/2}d\lambda\right\},
\eeq
which represents an even, $2\pi$--periodic function, holomorphic in the domain
$\dot{{\cI}}$ of the complex $\theta$--plane. Since in the following we shall
prove that this function can be represented by the Legendre expansion (\ref{sette}),
then it can properly be denoted by $f(\theta)$.
\end{proof}

Formula (\ref{quarantaquattro}) allows us to compute the boundary values
$f_\pm(v)$ (defined by $f_\epsilon(v)=\lim_{u\rightarrow 0^+}f(\epsilon u+iv)$;\,
$\epsilon= \pm$; $v\geq 0$) on the semiaxis $\{\theta=iv,\,v\geq 0\}$ in terms of
the corresponding boundary values $\hat{f}_\pm(w)$
(with $\gamma_{iv}:\{\tau=iw,\,0\leq w\leq v\}$), provided $\hat{f}_\pm(w)$ satisfy
a $C^1$--type regularity condition; the latter is definitely necessary in order
to perform the inversion of the Radon--Abel transform at the boundary.
The $C^1$--continuity of the boundary values follows from the fact that the
sequence $f_n=n^p a_n$, $p\geq 2$, satisfies the Hausdorff condition (\ref{tre}).
We thus obtain:
\beq
\label{quarantasette}
i[f_+(v)-f_-(v)]=F(v)=\frac{1}{\pi\sinh v}\frac{d}{dv}\int_0^v
e^{w/2}\,\hat{F}\,(w)\frac{\sinh w}{[2(\cosh v - \cosh w)]^{1/2}}\,dw,
\eeq
($\hat{F}(w)=i[\hat{f}_+(w)-\hat{f}_-(w)];\,\hat{f}_\epsilon(w)=
\lim_{t\rightarrow 0^+}\hat{f}(\epsilon t+iw); \,\epsilon= \pm$). \\
We can then apply the inverse Radon--Abel transform operator (defined by
formula (\ref{quarantaquattro})) to the series at the r.h.s. of formula
(\ref{ventisei}): i.e.,
\beq
\label{quarantotto}
\hat{f}(t)=\frac{1}{2\pi}e^{i(t-\pi)/2}
\sum_{n=-\infty}^{+\infty}(-1)^n\,a_n\,\cos\left[\left(n+\frac{1}{2}\right) (t-\pi)\right],
\eeq
and integrate term by term in view of the uniform convergence of this series,
which follows from the Hausdorff conditions on the coefficients $a_n$.
Next, we introduce the functions
\beq
\label{quarantottoprimo}
\begin{split}
\psi_n(\cos u) &=-\frac{i}{\pi\sin u}\frac{d}{du}\int_0^u
\cos\left[\left(n+\frac{1}{2}\right)(t-\pi)\right] \\
&\qquad\times\frac{\sin t}{\left[2(\cos u-\cos t)\right]^{1/2}}\,dt \qquad (0<u<2\pi),
\end{split}
\eeq
which are related to the Legendre polynomials $P_n(\cos u)$ as follows
(see formulae (II.79) and (II.91) of \cite{Bros}II):
\beq
\label{quarantanove}
\psi_n(\cos u)=\frac{(-1)^n}{4}(2n+1)P_n(\cos u).
\eeq
Finally, recalling that $a_n=-a_{-n-1}$ ($n\in\Z$), we obtain again the original
Legendre expansion (\ref{sette}):
\beq
\label{cinquanta}
f(u)=\uf(\cos u)=\frac{1}{\pi}\sum_{n=0}^\infty(-1)^n a_n\,\psi_n(\cos u)=
\frac{1}{4\pi}\sum_{n=0}^\infty(2n+1) a_n\,P_n(\cos u).
\eeq
We can restate the results of Proposition \ref{pro6} in the more natural geometry
of the $\cos\theta$--plane.

\begin{proposition}
\label{pro6primo}
If the sequence $f_n=n^p a_n$ $(n=0,1,\ldots)$ satisfies the Hausdorff condition
$(\ref{tre})$ with $p\geq 2$, then:
\begin{itemize}
\item[i)] the series $(\ref{I1})$ converges uniformly to an analytic function
$\uf(\cos u)$ $(u\equiv\Real\theta)$ in any compact domain $|\cos u| \leq |\cos u_0| < 1$; and
\item[ii)] the function $\uf(\cos u)$ admits a holomorphic extension to the
complex $\cos\theta$--plane $(\theta=u+iv)$ cut along the axis $[1,+\infty)$.
\end{itemize}
\end{proposition}

\section{Laplace transformation, representation of the jump function
and its reconstruction by the use of the Pollaczek polynomials}
\label{se:laplace}
From formulae (\ref{trenta}) and  (\ref{a.7}) the following equality follows:
\beq
\label{cinquantuno}
\tilde{a}(\lambda)=\int_0^{+\infty}e^{-(\lambda+1/2)w} (\cA F)(w)\,dw \qquad
\left (\lambda=\sigma+i\nu,\,\Real\lambda>-\frac{1}{2}\right).
\eeq
Writing explicitly the Abel transform $(\cA F)(w)$ (see formula (\ref{a.7})) yields
\beq
\label{cinquantadue}
\tilde{a}(\lambda)=2\int_0^{+\infty}e^{-(\lambda+1/2)w}
\left\{\int_0^w \frac{\underline{F}(\cosh v)\sinh v}{[2(\cosh w-\cosh v)]^{1/2}}\,dv\right\}\,dw
\quad \left(\Real\lambda>-\frac{1}{2}\right).
\eeq
If we exchange the integration order, this becomes
\beq
\label{cinquantatre}
\begin{split}
\tilde{a}(\lambda)=2\int_0^{+\infty} \underline{F}(\cosh v)\sinh v
\left\{\int_v^{+\infty} \frac{e^{-(\lambda+1/2)w}}{[2(\cosh w-\cosh v)]^{1/2}}\,dw\right\}\,dv \\
\qquad\left(\Real\lambda>-\frac{1}{2}\right).
\end{split}
\eeq
Recalling the integral representation of the second--kind Legendre functions
$Q_\lambda(\cosh v)$, i.e.,
\beq
\label{cinquantaquattro}
Q_\lambda(\cosh v)=\int_v^{+\infty} \frac{e^{-(\lambda+1/2)w}}{[2(\cosh w-\cosh v)]^{1/2}}\,dw
\qquad (\Real\lambda>-1, v>0),
\eeq
we can write formula (\ref{cinquantatre}) as follows
\beq
\label{cinquantacinque}
\tilde{a}(\lambda)=2\int_0^{+\infty}\underline{F}(\cosh v)Q_\lambda(\cosh v)\sinh v\,dv
\qquad\left(\Real\lambda>-\frac{1}{2}\right).
\eeq

\begin{remark}
The second--kind Legendre function presents a logarithmic singularity
at $v=0$; then the integral representation (\ref{cinquantaquattro}) holds true if
$v>0$; nevertheless, the integral in (\ref{cinquantacinque}) converges if
$\underline{F}(\cosh v)$ is regular at $v=0$.
\end{remark}

If $\Real\lambda=-1/2$ we can split $Q_{-1/2+i\nu}(\cosh v)$ into two terms,
$Q^{(E)}_{-1/2+i\nu}(\cosh v)$ and $Q^{(O)}_{-1/2+i\nu}(\cosh v)$, defined as follows:
\begin{eqnarray}
\label{cinquantaseia}
Q^{(E)}_{-1/2+i\nu}(\cosh v)&=&
\int_v^{+\infty}\frac{\cos\nu w}{[2(\cosh w-\cosh v)]^{1/2}}\,dw \qquad (v>0), \\
\label{cinquantaseib}
Q^{(O)}_{-1/2+i\nu}(\cosh v)&=&
-i\int_v^{+\infty}\frac{\sin\nu w}{[2(\cosh w-\cosh v)]^{1/2}}\,dw \qquad (v>0).
\end{eqnarray}
Next, we recall the following equality (see \cite{Bateman}):
\beq
\label{cinquantasette}
P_\lambda(\cos\theta)=\tan(\pi\lambda)\{Q_\lambda(\cos\theta)-Q_{-\lambda-1}(\cos\theta)\},
\eeq
(where we use a non--standard normalization of the $Q_\lambda$ functions,
which is more appropriate to our joint consideration of $P_\lambda$ and $Q_\lambda$;
the discrepancy with the usual notation is a factor $1/\pi$). We thus have the following
equality, which will be useful later on:
\beq
\label{cinquantotto}
P_{-1/2+i\nu}(\cosh v)=P_{-1/2-i\nu}(\cosh v)=2\tan\left[\pi\left(-\frac{1}{2}+i\nu\right)\right]
Q^{(O)}_{-1/2+i\nu}(\cosh v).
\eeq
We can now prove the following proposition.

\begin{proposition}
\label{pro7}
If the sequence $f_n=n^p a_n$ ($a_n$ being the Legendre coefficients)
satisfies the Hausdorff condition $(\ref{tre})$ with $p\geq 2$, then the jump
function $F(v)=\underline{F}(\cosh v)$ admits the integral representation
\beq
\label{cinquantanove}
F(v)=\underline{F}(\cosh v)=\frac{1}{4\pi}\int_{-\infty}^{+\infty}\tilde{a}(\sigma+i\nu)h(\sigma+i\nu)
P_{\sigma+i\nu}(\cosh v)\,d\nu \qquad \left(\sigma\geq-\frac{1}{2}\right),
\eeq
where $h(\sigma+i\nu)=2(\sigma+i\nu)+1$, and $P_{\sigma+i\nu}(\cosh v)$ denotes the
first--kind Legendre functions.
\end{proposition}

\begin{proof}
In Propositions \ref{pro3} and \ref{pro4} we derived the following formula:
\beq
\label{sessanta}
\hat{F}(w)=\frac{1}{2\pi}\int_{-\infty}^{+\infty}\tilde{a}(\sigma+i\nu)\,e^{(\sigma+i\nu)w}\,d\nu
\qquad \left(\sigma\geq-\frac{1}{2}\right).
\eeq
This formula can indeed be obtained by evaluating the discontinuity across $\tau=iw$,
$w\geq0$ of the function $g(\tau)-\exp(i\tau)g(-\tau)$, where
$g(\tau)=\frac{1}{2\pi}\sum_{n=0}^\infty a_n \exp(-in\tau)$. Note, in fact, that
the jump of $g(-\tau)$ across this half--line is equal to zero; but its expression
through a vanishing Cauchy integral allows us to introduce another equivalent integral
representation of $\hat{F}(w)$, i.e.,
\beq
\label{sessantuno}
\hat{F}(w)=\frac{i}{\pi}e^{-w/2}\int_{-\infty}^{+\infty}\tilde{a}(\sigma+i\nu)
\sin\left\{\left[\nu-i\left(\sigma+\frac{1}{2}\right)\right]w\right\}\,d\nu
\qquad \left(\sigma\geq-\frac{1}{2}\right).
\eeq
Applying the inverse Radon--Abel transform to $\exp(w/2)\hat{F}(w)$ (see formula
(\ref{quarantasette})) yields
\beq
\label{sessantadue}
\begin{split}
F(v)=-\frac{i}{\pi^2}\int_{-\infty}^{+\infty}\tilde{a}(\sigma+i\nu)\left\{
\frac{1}{\sinh v}\frac{d}{dv}\int_0^v\frac{\sin\left\{
\left[i\left(\sigma+\frac{1}{2}\right)-\nu\right]w\right\} \sinh w}
{\left[2(\cosh v-\cosh w)\right]^{1/2}}\,dw\right\}\,d\nu \\
\qquad\left(\sigma\geq-\frac{1}{2}\right).
\end{split}
\eeq
The r.h.s. of formula (\ref{sessantadue}) converges to $F(v)$ if $\hat{F}(w)$ is
of class $C^1$ and $\nu\tilde{a}(\sigma+i\nu)$ belongs to $L^1(-\infty,+\infty)$.
Both properties follow from the requirement that the sequence $f_n=n^p a_n$ satisfies
the Hausdorff condition (\ref{tre}) with $p\geq 2$. Finally we recognize in the
integrand of formula (\ref{sessantadue}) the first--kind Legendre function
$P_{\sigma+i\nu}(\cosh v)$; in fact, we have (see formula (II.86) of \cite{Bros}II):
\beq
\label{sessantatre}
\begin{split}
&\frac{1}{4\pi}P_{\sigma+i\nu}(\cosh v)[2(\sigma+i\nu)+1] \\
&\qquad =-\frac{i}{\pi^2\sinh v}\frac{d}{dv}
\int_0^v \frac{\sin\left\{\left[i\left(\sigma+\frac{1}{2}\right)-\nu\right]w\right\}}
{\left[2(\cosh v-\cosh w)\right]^{1/2}}\sinh w\,dw.
\end{split}
\eeq
By plugging (\ref{sessantatre}) into (\ref{sessantadue}), we get the result,
i.e., formula (\ref{cinquantanove}).
\end{proof}

In the particular case of $\sigma=-1/2$, in view of the evenness of the function
$P_{-1/2+i\nu}(\cosh v)$ with respect to $\nu$, only the odd component of
$\tilde{a}(-1/2+i\nu)$ contributes to the integral in (\ref{cinquantanove}).
Accordingly, in view of formula (\ref{cinquantotto}), we can write the Laplace
transform (\ref{cinquantacinque}) in terms of the function
$P_{-1/2+i\nu}(\cosh v) / \tan[\pi(-1/2+i\nu)]$, instead of $Q_{-1/2+i\nu}(\cosh v)$.
It can easily be verified that, in this case, formulae (\ref{cinquantacinque})
and (\ref{cinquantanove}) give (up to normalization constants) the classical
Mehler transform (see \cite{Bateman}), which is precisely the tool used by Stein and
Wainger \cite{Stein} for proving their theorem.

We can now rapidly mention how the discontinuity function can be reconstructed,
starting from the Fourier coefficients, by the use of the Pollaczek polynomials.

\begin{proposition}
\label{pro8}
Let us suppose that the sequence of the Fourier--Legendre coefficients $a_n$,
$(n=0,1,2,\ldots)$ satisfies the Hausdorff condition $(\ref{tre})$; then the
function $\hat{F}(w)e^{w/2}$ (see formula $(\ref{trentuno})$), can be
represented by the following expansion which converges in the sense of the
$L^2$--norm:
\beq
\label{la1}
\hat{F}(w)e^{w/2} = \sum_{\ell=0}^\infty c_\ell \Phi_\ell(w) \qquad (w\in\R^+),
\eeq
where
\beq
\label{la2}
c_\ell =
\sqrt{2}\sum_{n=0}^\infty\frac{(-1)^n}{n!}a_n \cP_\ell\left[-i\left(n+\frac{1}{2}\right)\right],
\eeq
\beq
\label{la3}
\Phi_\ell(w)=i^\ell\sqrt{2}e^{-w/2}L_\ell(2e^{-w})e^{-e^{-w}},
\eeq
where $\cP_\ell$ and $L_\ell$ are the Pollaczek and the Laguerre polynomials
\cite{Bateman}, respectively.
\end{proposition}

\begin{proof}
Let us note that
$\tilde{a}(-1/2+i\nu)$ and $\hat{F}(w)e^{w/2}$ belong to $L^2(-\infty,+\infty)$
(see statement (iii) of Proposition 1 and formula (\ref{trentuno})). We can then write
\beq
\label{la4}
\hat{F}(w)e^{w/2}=  \lm\displaylimits_{\nu_0\rightarrow +\infty}
\left(\frac{1}{2\pi}\int_{-\nu_0}^{\nu_0}\tilde{a}\left(-\frac{1}{2}+i\nu\right) e^{i\nu w}\,d\nu\right),
\eeq
and the proof proceeds exactly as in Theorem 2 of \cite{DeMicheli1}, where
$\hat{F}(w)e^{w/2}$ now plays the role of $F(v)e^{v/2}$ (see also the Remark at
the end of Proposition 1).
\end{proof}

\begin{remark}
The previous result holds true under the condition that the sequence
$f_n=n^p a_n$, ($n=0,1,\ldots$) satisfies condition (\ref{tre}) with $p\geq 0$;
but in the next proposition, which involves the Abel transform and its inverse,
the latter condition must be satisfied with $p\geq 2$ (see formulae
(\ref{quarantasette}) and (\ref{a.7})).
\end{remark}

Next we set (see formula (\ref{a.7})): $\cosh v = y$, $\cosh w = x$; we can rewrite
the Abel transform as a convolution product of the following form
\beq
\label{la5}
(\cA F)(x)=\sqrt{2}\int_1^x \frac{\underline{F}(y)}{\sqrt{(x-y)}}\,dy :\equiv \sqrt{2}\,(\underline{F} * \chi_1)(x).
\eeq
From formula (\ref{quarantasette}) we analogously obtain
\beq
\label{la6}
\underline{F}(y) = \frac{1}{\pi\sqrt{2}}\frac{d}{dy}\int_1^y \frac{(\cA F)(x)}{\sqrt{(y-x)}}\,dx :\equiv
\frac{1}{\pi\sqrt{2}}\frac{d}{dy}(\cA F * \chi_1)(y).
\eeq
We can then prove the following proposition.

\begin{proposition}
\label{pro9}
If the sequence $f_n=n^p a_n$ $(n=0,1,2,\ldots)$ satisfies the Hausdorff
condition $(\ref{tre})$ with $p\geq 2$, then the following limit holds true:
\beq
\label{la7}
\lim_{m\rightarrow\infty} \langle(\underline{F}-\psi_m),\phi\rangle = 0,
\eeq
where
\beq
\label{la8}
\psi_m=\frac{1}{\pi\sqrt{2}}\frac{d}{dy}\left[\left(\sum_{\ell=0}^m
c_\ell\Phi_\ell\right) * \chi_1\right],
\eeq
$\phi\in S_\infty(\R)$ ($S_\infty(\R)$ being the Schwartz space of the $C^\infty(\R)$
functions $\phi(x)$ that together with all their derivatives decrease, for $|x|$
tending to $\infty$, faster than any negative power of $|x|$), and $\langle f,\phi\rangle$
denotes the Lebesgue integral $\int_{-\infty}^{+\infty} \bar{f}\phi\,dx$.
\end{proposition}

\begin{proof}
We have from formulae (\ref{la6}) and (\ref{la8}):
\beq
\label{la9}
\begin{split}
\left\langle(\underline{F}-\psi_m),\phi\right\rangle &=
\frac{1}{\pi\sqrt{2}}\left\langle\frac{d}{dy}\left[\left(\cA F-\sum_{\ell=0}^m
c_\ell\Phi_\ell\right)* \chi_1 \right],\phi\right\rangle \\
&=-\frac{1}{\pi\sqrt{2}}\left\langle\left(\cA F-\sum_{\ell=0}^m
c_\ell\Phi_\ell\right)* \chi_1,\phi'\right\rangle \qquad (\phi\in S_\infty(\R)).
\end{split}
\eeq
From inequality (\ref{ventotto}) and formula (\ref{a.7}) it follows that $\cA F$
has a power--like behavior in $x$.
Next, in view of the Fubini theorem, we have:
\beq
\label{la10}
\left\langle\left(\cA F-\sum_{\ell=0}^m c_\ell\Phi_\ell\right)* \chi_1,\phi'\right\rangle=
\left\langle\cA F-\sum_{\ell=0}^m c_\ell\Phi_\ell,\phi' * \chi^\infty\right\rangle,
\eeq
where
\beq
\label{la11}
(\phi' * \chi^\infty)(x)=\int_x^{+\infty}\frac{\phi'(y)}{\sqrt{(y-x)}}\,dy < \infty.
\eeq
From the Schwarz inequality it follows:
\beq
\label{la12}
\left\langle\cA F-\sum_{\ell=0}^m c_\ell\Phi_\ell,\phi' * \chi^\infty\right\rangle \leq
\left\|\cA F-\sum_{\ell=0}^m c_\ell\Phi_\ell\right\|_{L^2[0,+\infty)} \cdot
\left\|\phi' * \chi^\infty\right\|_{L^2[0,+\infty)}.
\eeq
Now, in view of the fact that
$\lim_{m\rightarrow\infty}\left\|\cA F-\sum_{\ell=0}^m c_\ell\Phi_\ell\right\|_{L^2[0,+\infty)}=0$
(see Proposition \ref{pro8}) and that
$\left\|\phi' * \chi^\infty\right\|_{L^2[0,+\infty)}<\infty$,
statement (\ref{la7}) holds true.
\end{proof}

\appendix
\section{Horocycles and Radon--Abel transformations
on the (real and complexified) one--sheeted hyperboloids}
Let $X$ denote the real one--sheteed hyperboloid in $\R^3$, with equation:
\beq
\label{a.1}
x_0^2-x_1^2-x_2^2=-1 \qquad (x \equiv (x_0,x_1,x_2)),
\eeq
and let $\hat{X}$ be the meridian hyperbola, in the $x_1=0$ plane, with
equation $x_0^2-x_2^2=-1$ (see Fig. \ref{figura_1}).
The manifold $X$ can be described by the following polar coordinates:
\begin{subequations}
\label{a.2}
\begin{eqnarray}
x_0 &=& \sinh v \cosh\phi \qquad (v,\phi\in\R), \label{a.2a} \\
x_1 &=& \sinh v \sinh\phi, \label{a.2b} \\
x_2 &=& \cosh v. \label{a.2c}
\end{eqnarray}
\end{subequations}

\begin{figure}[tb]
\begin{center}
\includegraphics[width=9cm]{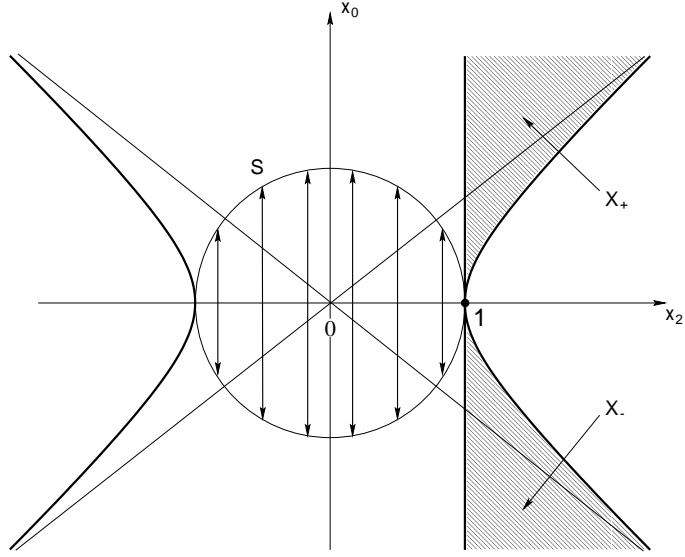}
\caption{\label{figura_2}$x_1=0$ section of the one--sheteed hyperboloid.}
\end{center}
\end{figure}

Hereafter we shall deal more specifically with the following region of $X$
(see Fig. \ref{figura_2}):
\beq
\label{a.3}
X_+=\{x\in X\,:\,x_0\geq 0,\,x_2\geq 1\}.
\eeq
Besides polar coordinates, another system of local coordinates on $X$ is equally
valid for describing the set $X_+$, namely the horocyclic coordinates:
\begin{subequations}
\label{a.4}
\begin{eqnarray}
x_0 &=& \sinh w +\frac{1}{2}\zeta^2\,e^w \qquad (\zeta\in\R,\,w\in\R^+), \label{a.4a} \\
x_1 &=& \zeta\,e^w, \label{a.4b} \\
x_2 &=& \cosh w - \frac{1}{2}\zeta^2\,e^w. \label{a.4c}
\end{eqnarray}
\end{subequations}
The sections $w = \Cons$ are parabolae lying in the planes $x_0+x_2=e^w$, called
horocycles.

Next we introduce the following integral:
\beq
\label{a.5}
\int_{h_w}\underline{F}\left(\cosh w-\frac{1}{2}\zeta^2\,e^w\right)\,d\zeta = \hat{F}(w),
\eeq
where $h_w$ is the oriented segment of horocycle belonging to $X_+$, which is
represented  by the arc of the parabola whose apex, which lies on $\hat{X}$, is
obtained by setting $\zeta=0$ in eqs. (\ref{a.4}) (i.e., with coordinates
$x_0=\sinh w,\,x_1=0,\,x_2=\cosh w$), and whose endpoints lie on the plane $x_2=1$.
Moreover, the function $\underline{F}$ is assumed to satisfy the regularity conditions
that make the integral (\ref{a.5}) convergent. Since the integrand is an even
function of $\zeta$, the integration domain can be restricted to the part of $h_w$
with $x_1\geq 0$, which is given by
\beq
\label{a.6}
h_w^+:\left\{\zeta=\left[2e^{-w}\,(1-\lambda)(\cosh w-1)\right]^{1/2},
\,(0\leq\lambda\leq 1),\,w\in\R^+\right\}.
\eeq
For $\lambda=1$ we get $\zeta=0$, i.e., the apex of the parabola representing
the horocycle; for $\lambda=0$ we get $\zeta=\left[2e^{-w}(\cosh w-1)\right]^{1/2}$,
which gives the intersection of the horocycle with the plane $x_2=1$.

Integral (\ref{a.5}) defines a transformation of Radon--type in $X$, where the
horocycles play the same role as the planes do in the ordinary Radon transformation.
Next, if we set $x_2=\cosh v$, we have
$\zeta(v)=\left[2e^{-w}(\cosh w-\cosh v)\right]^{1/2}$, which varies again between
$\zeta=0$ for $v=w$ (apex of the parabola) and $\zeta=\left[2e^{-w}(\cosh w-1)\right]^{1/2}$
for $v=0$ (endpoint of the parabola: i.e., $x_2=\cosh v=1$). Since
$d\zeta/dv=-e^{-w/2}\sinh v [2(\cosh w-\cosh v)]^{-1/2}$, integral (\ref{a.5})
takes the form:
\beq
\label{a.7}
\begin{split}
\hat{F}(w) &= 2e^{-w/2}\int_0^w \underline{F}(\cosh v)\frac{\sinh v}{[2(\cosh w-\cosh v)]^{1/2}}\,dv \\
& :\equiv e^{-w/2} (\cA F)(w)~,~~~~(w\in\R^+),
\end{split}
\eeq
which is an Abel--type integral.

We notice that the set of horocycles $\{h_w;\,w\geq 0\}$ defines a fibration with
basis $\hat{X}$ on the domain $X$. We denote by $\underline{h}$ the projection
associated with this fibration; in particular, all the points of the domain $X_+$
are projected on the basis $\hat{X}$, i.e., $\forall x \in X_+,(x=(x_0,x_1,x_2)),
\underline{h}(x)=x_w$ is the intersection of $\hat{X}$ with the unique horocycle
$h_w$ which contains $x$.

We can now regard the one--sheeted hyperboloid $X$ as a real submanifold of a
complex hyperboloid $X^{(c)}$ (see Fig. \ref{figura_1}), whose equation is given by:
\beq
\label{a.8}
z_0^2-z_1^2-z_2^2=-1 \qquad (z=(z_0,z_1,z_2),\,z_i\in\C,\,i=0,1,2).
\eeq
Accordingly, we introduce the following complex--valued polar coordinates:
\begin{subequations}
\label{a.9}
\begin{eqnarray}
z_0 &=& -i\sin\theta\cosh\phi \qquad (\theta,\phi\in\C), \label{a.9a} \\
z_1 &=& -i\sin\theta\sinh\phi, \label{a.9b} \\
z_2 &=& \cos\theta. \label{a.9c}
\end{eqnarray}
\end{subequations}
If we set $\theta=iv$ ($v\in\R$) in eqs. (\ref{a.9}), and assume $\phi$ real, we
obtain the polar coordinates (\ref{a.2}) that describe the real hyperboloid $X$;
if we set $\theta=u$ ($u\in\R$) and $\phi=i\eta$ ($\eta\in\R$), we obtain the
\emph{Euclidean} sphere:
\begin{subequations}
\label{a.10}
\begin{eqnarray}
z_0 &=& -i\sin u \cos\eta \qquad (u,\eta\in\R), \label{a.10a} \\
z_1 &=& \sin u\sin\eta, \label{a.10b} \\
z_2 &=& \cos u. \label{a.10c}
\end{eqnarray}
\end{subequations}
We now want to extend the fibration with basis $\hat{X}$, introduced above, to the
complex one--sheeted hyperboloid $X^{(c)}$. For this purpose, we consider the
complex meridian hyperbola $\hat{X}^{(c)}$, lying in the $z_1=0$ plane, with
equation: $z_0^2-z_2^2=-1$. The intersections of $\hat{X}^{(c)}$ with the family
of planes $P_\tau$, with equations $z_0+z_2=e^{-i\tau}$ ($\tau\in\C,\,\tau=t+iw$)
are the points $z_\tau$ whose coordinates are $(z_\tau)_0=-i\sin\tau$, $(z_\tau)_2=\cos\tau$.
The sections of $X^{(c)}$ by the planes $P_\tau$ are complex parabolae (except in
the case $z_0+z_2=0$) lying in the planes $P_\tau$. These complex parabolae are the
geometric realization of the complex horocycles, which will be denoted by $h_\tau$ $(\tau\in\C)$.
At this point, it is convenient to introduce the domain
$X'^{(c)}=\{z\in X^{(c)};\,z_0+z_2 \neq 0\}$ (dense in $X^{(c)}$),
which can be parametrized by the horocyclic coordinates $(\zeta,\tau)$ in the following way:
\begin{subequations}
\label{a.11}
\begin{eqnarray}
z_0 &=& -i\sin\tau + \frac{1}{2}\zeta^2e^{-i\tau} \qquad (\zeta,\tau\in\C), \label{a.11a} \\
z_1 &=& \zeta e^{-i\tau}, \label{a.11b} \\
z_2 &=& \cos\tau - \frac{1}{2}\zeta^2e^{-i\tau}. \label{a.11c}
\end{eqnarray}
\end{subequations}
We can now introduce the cut--domain $X^{(c)}\setminus\Sigma^{(c)}$,
where $\Sigma^{(c)}$ is defined as follows:
\beq
\label{a.12}
\Sigma^{(c)}=\{z\in X^{(c)};\, z_2\in [1,+\infty)\}.
\eeq
The domain $\Sigma^{(c)}\cap X$ is composed by two sets, i.e. $X_+$ (defined by
(\ref{a.3})) and $X_-$ defined by (see Fig. \ref{figura_2}):
\beq
\label{a.13}
X_-=\{x\in X;\,x_0\leq 0,\,x_2\geq 1\}.
\eeq
The fibration produced by the horocycles $h_w$ can now be extended through the
use of complex horocycles $h_\tau$, whose intersections with the meridian
(complex) hyperbola $\hat{X}^{(c)}$ are the points $z_\tau$ with coordinates
$(z_\tau)_0=-i\sin\tau,\,(z_\tau)_2=\cos\tau$.
In the following we shall deal with functions that only depend on the coordinate
$z_2=\cos\theta$, i.e., $\uf(\cos\theta)$ (or, alternatively, $f(\theta)$).
Furthermore, we assume that these even and $2\pi$--periodic functions $f(\theta)$
are holomorphic in the cut--domain
$\dot{{\cI}}=\cI^{(0)}_+\setminus\dot{\Xi}^{(0)}_+ \cup\cI^{(0)}_-\setminus\dot{\Xi}^{(0)}_-$,
and continuous up to the boundaries of this domain (notice that
the notations $\cI^{(0)}_\pm$ and $\dot{\Xi}^{(0)}_\pm$ have been introduced in
connection with Proposition \ref{pro3} where they referred to the complex plane
of the variable $\tau$). Moreover, we denote by $\underline{D}\in\C$ the domain
of analyticity of $\uf(\cos\theta)$ in the $\cos\theta$--plane. Next, we introduce
the following integral:
\beq
\label{a.14}
2\int_{h^+_\tau} \uf\left(\cos\tau-\frac{1}{2}\zeta^2\,e^{-i\tau}\right)\,d\zeta = \hat{f}(\tau),
\eeq
$h^+_\tau$ being the arc of complex horocycle defined by
\beq
\label{a.15}
h^+_\tau: \left\{\zeta=\left[2\,e^{i\tau}(1-\lambda)(\cos\tau-1)\right]^{1/2},\,
(0\leq\lambda\leq 1),\, \tau\in\C\right\}.
\eeq
For $\lambda=1$, we get $\zeta=0$: i.e., the point $z_\tau$ (belonging $\hat{X}^{(c)}$);
for $\lambda=0$, we get $\zeta=[2\exp(i\tau)(\cos\tau-1)]^{1/2}$, which gives the
intersection of $h^+_\tau$ with the plane $z_2=1$ (in particular, if we set $\tau=iw$,
$w\in\R^+$, we obtain the expression $\zeta(w)=[2\exp(-w)(\cosh w-1)]^{1/2}$,
previously established in connection with formula (\ref{a.6})).
Then, by setting $z_2=\cos\theta$ $(\theta\in\C,\,\theta=u+iv)$, we have:
$\zeta(\theta)=[2e^{i\tau}(\cos\tau-\cos\theta)]^{1/2}$, and
$d\zeta=e^{i\tau/2}[2(\cos\tau-\cos\theta)]^{-1/2}\sin\theta\,d\theta$.
Furthermore, since on $h^+_\tau$ we have:
$\cos\theta(\lambda)-1=\lambda(\cos\tau-1)$ $(0\leq\lambda\leq 1)$, integral (\ref{a.14})
can be rewritten in the following form
\beq
\label{a.16}
\hat{f}(\tau)=
-2e^{i\tau/2}\int_{\gamma_\tau}f(\theta)\left[2(\cos\tau-\cos\theta)\right]^{-1/2}\sin\theta\,d\theta,
\eeq
where $\gamma_\tau$ denotes the \emph{ray} $\underline{\gamma}_\tau$ oriented
from 0 to $\tau$, and
$\underline{\gamma}_\tau:\{\theta=\theta(\lambda);\,\cos\theta(\lambda)-1=\lambda(\cos\tau-1),\,
0\leq\lambda\leq 1,\,\theta(0)=0,\,\theta(1)=\tau\}$.
Moreover, the relevant branch of the function $[2(\cos\tau-\cos\theta)]^{-1/2}$ is
specified by the condition that for $\tau=iw$, and $\theta=iv$ (with $w>v$), it takes
the value $[2(\cosh w-\cosh v)]^{-1/2}\geq 0$. In fact, when $\tau=iw$ ($w>0$),
the horocycle $h_\tau=h_{iw}$ is real and carried by the hyperboloid $X$.
Moreover, in this case, transformation (\ref{a.16}) can be applied to the boundary
values of $f$ on the opposite sides of the cut (corresponding to the domain $X_+$),
and, in particular, to the corresponding discontinuity function; in this way formula
(\ref{a.7}) is reobtained, and the function $\underline{F}$ now represents the jump
function across the cut.

We can now show that if $f(\theta)$ is an even $2\pi$--periodic function
holomorphic in the cut--domain $\dot{{\cI}}$, then $\hat{f}(\tau)$ is a
$2\pi$--periodic function holomorphic in $\dot{{\cI}}$ that satisfies the
following symmetry relation:
\beq
\label{a.17}
\hat{f}(\tau)=-e^{i\tau}\,\hat{f}(-\tau).
\eeq
In order to prove this statement we rewrite expression (\ref{a.16}) in the
following form
\beq
\label{a.18}
\hat{f}(\tau)=e^{i\tau/2}\left[2(\cos\tau-1)\right]^{1/2}\int_0^1
\uf(1+\lambda(\cos\tau-1))(1-\lambda)^{-1/2}\,d\lambda,
\eeq
where the following parametrization $\cos\theta(\lambda)=1+\lambda(\cos\tau-1)$
in integral (\ref{a.16}) has been used. Since the integral in (\ref{a.18}) is a
function of $\cos\tau$ analytic in $\underline{D}$, $\hat{f}(\tau)$ can be written as
\beq
\label{a.19}
\hat{f}(\tau)=e^{i\tau/2}\left(\sin\frac{\tau}{2}\right)\,a(\cos\tau),
\eeq
where $a(\cos\tau)$ is a function holomorphic in $\underline{D}$. From
representation (\ref{a.19}) one recovers that $\hat{f}(\tau)$ is $2\pi$--periodic
and, therefore, holomorphic in $\dot{{\cI}}$, and, in addition, the symmetry
relation (\ref{a.17}) is satisfied.

Finally, by restricting formulae (\ref{a.14}) and (\ref{a.16}) to the set of
real values of the variables $\tau$ and $\theta$, namely $\tau=t$, $\theta=u$,
from (\ref{a.16}) we obtain:
\beq
\label{a.20}
\hat{f}(t)=-2\,e^{it/2}\int_0^t\,f(u)\left[2(\cos t-\cos u)\right]^{-1/2}\sin u\,du.
\eeq
By taking into account the relevant branch of the factor $[2(\cos t-\cos u)]^{-1/2}$,
formula (\ref{a.20}) can be written in the following more precise form (involving
a positive bracket),
\beq
\label{a.21}
\hat{f}(t)=-2i\epsilon(t)\,e^{it/2}\int_0^t\,f(u)\left[2(\cos u-\cos t)\right]^{-1/2}\sin u\,du,
\eeq
where $\epsilon(t)$ denotes the sign function. Note that formula (\ref{a.21})
coincides with formula (\ref{undici}).

\subsection*{Acknowledgments}
The authors feel deeply indebted to Prof. J. Bros for many illuminating discussions.

\end{document}